\documentclass{amsart}

\input{xy}
\xyoption{all}
\newcommand{\ot}{\leftarrow}

\newcommand{\CA}{\mathcal{C}_A}
\newcommand{\DA}{\mathcal{D}^b(\textup{mod}\,A)}

\DeclareMathOperator{\Hom}{Hom}%
\DeclareMathOperator{\Ext}{Ext}%
\DeclareMathOperator{\ind}{ind}%
\DeclareMathOperator{\End}{End}%

\newcommand{\za}{\alpha}
\newcommand{\zb}{\beta}
\newcommand{\zd}{\delta}

\newcommand{\zg}{\gamma}

\newcommand{\zs}{\sigma}

\newtheorem{thm}{Theorem}[section]

\newtheorem{cor}[thm]{Corollary}
\newtheorem{lem}[thm]{Lemma}

\newtheorem{definition}{Definition} 
\newtheorem{rem}[thm]{Remark}

\newenvironment{pf}{{Proof}.}

\title{Exceptional sequences and clusters}
\author{Kiyoshi Igusa and Ralf Schiffler}
\thanks{The second author is supported by the NSF grant DMS-0908765 and the University of Connecticut}

\subjclass[2000]{
16G20; 20F55}

\keywords{
cluster category, cluster-tilting set, braid group, Dynkin quiver, reflection groups, real Schur roots
}

\begin{document}

\begin{abstract}
We show that exceptional sequences for hereditary algebras are characterized by the fact that the product of the corresponding reflections is the inverse Coxeter element in the Weyl group. We use this result to give a new combinatorial characterization of clusters tilting sets in the cluster category in the case where the hereditary algebra is of finite type.
\end{abstract}

\maketitle




\section*{Introduction}


Let $W$ be a Coxeter group and let $S=\{s_1,s_2,\ldots,s_n\}$ be the set
of simple reflections. Denote by $C$ the Coxeter element $C=s_1s_2\cdots
s_n$. As a first main result of this paper, we show that \emph{the braid
group acts transitively on the set of all sequences of $n$ reflections
whose product is $C$}, see Theorem \ref{thm 2}.

We then use this result in the crystallographic case to investigate the
relation between exceptional sequences and clusters, first in the case of
a hereditary algebra of finite type and then over the path algebra of an
arbitrary quiver without oriented cycles. This is done as follows.

The sequence of simple modules in reverse order: $(S_n,\cdots,S_1)$ is a
complete exceptional sequence if the projective cover of each $S_j$
contains only $S_i$ for $i\le j$ in its composition series.
Crawley-Boevey \cite{C-B93} and Ringel \cite{R94} showed that the braid
group acts transitively on the set of complete exceptional sequences of
indecomposable modules over a hereditary algebra.
This action preserves the product of the corresponding reflections in the
Weyl group, thus,
for any complete exceptional sequence $(E_1,\cdots,E_n)$, the product of
the corresponding reflections is equal to the inverse Coxeter element
$C^{-1}$, since the latter is the product of reflections corresponding to
the exceptional sequence given by the simple modules in reverse order:
\[
	s_{E_1}s_{E_2}\cdots s_{E_n}=C^{-1}=(s_1s_2\cdots s_n)^{-1}=s_n\cdots s_1.
\]
By our Theorem \ref{thm 2}, it follows that this equation holds \emph{ if
and only if}  $(E_1,\cdots,E_n)$ is a complete exceptional sequence.

We note that these results are known in the finite
case \cite{Bessis03}, \cite{BW08} and have been extended to the
affine case by Ingalls and Thomas \cite{Ing-Th}. Our new proof is type
independent and includes these previous results as corollaries.
(We note however that, by \cite{Ing-Th},
the lattice condition proved in \cite{BW08} does
not hold in general.)

When $A$ is a hereditary algebra of finite type, the Weyl group is finite and has a unique element $w_0$ of maximal length. This is the element which sends all positive roots to negative roots and vice versa. It can be written as a product of simple reflections $s_i$, one for every indecomposable module in the $\tau$ orbit of the $i$th projective module $P_i$, where $\tau$ is Auslander-Reiten translation. When these simple reflections are written in adapted order (in the order they occur in the Auslander-Reiten quiver) then we get a reduced expression $s_{i_1}s_{i_2}\cdots s_{i_\nu}$ for $w_0$. The cluster category $\mathcal C_A$ of $A$ contains $n$ more indecomposable objects given by the shifted projective modules $P_i[1]$. When we add the corresponding simple reflections to this reduced expression for $w_0$ on the right we get an unreduced expression  $s_1s_2\cdots s_n s_{i_1}s_{i_2}\cdots s_{i_\nu}$ for the element $w_1=Cw_0\in W$.
Each indecomposable object in the cluster category then corresponds to exactly one simple reflection in this expression for $w_1$, for example $s_{i_1}$ corresponds to $\tau^{-1}P_1$. If $1\le t_1<\cdots<t_n\le n+\nu$, we denote by  $w^\delta(t_1,\cdots,t_n)$ the word obtained from  $s_1s_2\cdots s_n s_{i_1}s_{i_2}\cdots s_{i_\nu}$ by deleting the $n$ simple reflections at the positions $t_1,\ldots,t_n$. Call this the \emph{deleted word}. We show that  $w^\delta(t_1,\cdots,t_n)$ is a reduced expression for $w_0$ if and only if the corresponding set of $n$ indecomposable objects in the cluster category $\mathcal{C}_A$ is a cluster-tilting set, see Theorem \ref{thm main}.
We also describe the mutations in terms of the reduced expression  $w^\delta(t_1,\cdots,t_n)$, see Theorem \ref{thm 3}.

The paper is organized as follows. In section \ref{sect 1.1} we prove that the braid group acts transitively on the set of all sequences of $n$ reflections whose product is the inverse Coxeter element $C^{-1}$. The argument is a generalization of our deleted word construction to the infinite case. In section \ref{sect main}, we precisely formulate the statement that a set of $n$ objects of a cluster category of finite type in adapted order forms a cluster tilting set if and only if the corresponding deleted word is reduced. In section \ref{sect proof} we prove this by observing that the first condition is equivalent to the condition that the sequence of objects gives an exceptional sequence and the second condition is equivalent to the condition that the product of the corresponding sequence of reflections is $C^{-1}$. One key idea is that of ``algebraic mutation'' which parallels mutation of clusters. In subsection \ref{example} we illustrate our theorem in type $A_n$ with an example and in subsection \ref{Woo} we use the example to give a cluster-tilting theoretic interpretation of a result of Woo \cite{Woo04}.

In section \ref{sect infinite type} we derive the corollary that a real root is a real Schur root if and only if the corresponding reflection is a prefix of the Coxeter element.

At the end of the paper we include an appendix communicated to us by Hugh Thomas, in
which the main theorem of the paper is used to prove that the set of finitely generated, exact abelian, extension-closed subcategories of the module category of a path algebra over a quiver without oriented cycles is in bijection with the set of prefixes of the corresponding Coxeter element.


\begin{section}{Braid group actions in Coxeter groups}


\begin{subsection}{Coxeter groups}\label{sect 1.1}

Let $W$ be a Coxeter group, and let  $S=\{s_1,s_2,\ldots,s_n\}$ be the set of simple reflections. $W$ is generated by $S$ subject to the relations $(s_is_j)^{m_{ij}}=1$, for some $m_{ij}$ such that $m_{ii}=1$ and $m_{ij}\ge 2$ if $i\ne j$. (See \cite{H} or \cite{Bro} for basic properties of Coxeter groups including all that we will be using.)

We use the standard bilinear pairing $B$ on $\mathbb{R}^n$ given on the standard unit vectors $\za_i$ by $B(\za_i,\za_j)=-\cos(\pi/m_{ij})$ when $m_{ij}<\infty$ and $B(\za_i,\za_j)=-1$ if $m_{ij}=\infty$. It is well known, for example by \cite{H}, 5.3, that this pairing gives a faithful linear action of the Coxeter group $W$ on $\mathbb{R}^n$ by the formula:
\[
	s_i(x)=x-2B(\za_i,x)\za_i.
\]

The \emph{root system} $\Phi\subset\mathbb{R}^n$ is the set of all $w(\za_i)$ where $w\in W$ and $i=1,2,\ldots,n$. Every root is a positive or negative linear combinations of simple roots: $\Phi=\Phi_+\coprod\Phi_-$.

The set of \emph{reflections} $T\subseteq W$ is defined as 
\[ T=\bigcup_{w\in W} wSw^{-1} .\]
The Weyl group $W$ acts on $T$ by conjugation.

There is a bijection $\Phi_+\cong T$ given by sending $\za\in\Phi_+$ to $s_\za\in T$ given by
\[
	s_\za(x)=x-2B(\za,x)\za.
\]
Lemma \ref{lem 1} shows that $s_\za\in T$ and that this mapping $\Phi_+\to T$ is $W$-equivariant. The inverse mapping $T\to \Phi_+$ sends the reflection $t\in T$ to the unique unit vector $\za$ which is a nonnegative linear combination of the simple roots $\za_i$ so that $t(\za)=-\za$.

\begin{lem}\label{lem 1}
For any root $\za$ and any $w\in W$ we have
\[ws_\za w^{-1}=s_{w(\za)}.\]
\end{lem}

\begin{pf} \cite[5.7]{H} If $y=w(x)$ then
$
ws_\za w^{-1}(y)=w(x-2B(\za,x)\za)=w(x)-2B(\za,x)w(\za)=w(x)-2B(w(\za),w(x))w(\za)=s_{w(\za)}(y)
$ since $B$ is $W$-invariant.
\qed\end{pf}

Note that $s_{-\za}=s_{\za}$. We also note that the set of reflections depends only on the pair $(W,S)$ whereas the set of positive roots depends on our choice of linear action of $W$ on $\mathbb{R}^n$ and our arbitrary decision to make them all unit vector. If we modify our choice of positive roots by multiplying them with positive scalar, making sure that roots in the same orbit of the action of $W$ are multiplied by the same scalar and define $s_{r\za}=s_\za$ for $r\neq0$, none of the statements below will be affected.

Let $C$ be the \emph{Coxeter element} $C=s_1s_2\cdots s_n$. Since the numbering of the simple reflections is arbitrary, $C$ represents the product of the elements of $S$ in any fixed order. We know that $C$ and $C^{-1}$ have length $n$. So, there are exactly $n$ positive roots $p_1,p_2,\ldots,p_n$ which are sent to negative roots by $C^{-1}$.  We call these the \emph{projective roots}.

\end{subsection}

\begin{subsection}{Braid group action}

Let $B_m$ be the braid group on $m$ strands and denote its generators by $\zs_1,\zs_2,\ldots,\zs_{m-1}$. Then $B_m$ acts on the set of all $m$ element sequences in any group as follows: the generator $\zs_i$ acts by moving $g_i$ one space to the right and conjugating $g_{i+1}$ by $g_{i}$:
\[
	\zs_i(g_1,\cdots,g_m)=(g_1,\cdots,g_{i-1},g_ig_{i+1}g_i^{-1},g_i,g_{i+2},\cdots,g_m).
\]
Note that the product of the group elements remains the same. Also note that, for any conjugacy class $X$, the set of sequences $X^m$ is invariant under the action of $B_m$.
The braid group $B_m$ also acts on the set of sequences of $m$ positive roots by
\[
	\zs_i(\zb_1,\cdots,\zb_m)=(\zb_1,\cdots,\zb_{i-1},\zb_{i+1}',\zb_i,\zb_{i+2},\cdots,\zb_m),
\]
where $\beta'_{i+1}= \vert s_{\beta_i}(\beta_{i+1} )\vert \in\Phi_+$, here we use the notation $\vert\beta\vert = \beta $, if $\beta \in \Phi_+$, and $\vert\beta\vert = -\beta $, if $\beta \in \Phi_-$. Note that $
s_{\zb_{i+1}'}=s_{\zb_i}s_{\zb_{i+1}}s_{\zb_i}
$ and therefore the action of $B_m$ on $\Phi_+^m$ agrees with the action on $T^m$, that is, the bijection $\Phi_+^m\cong T^m$ sending $(\zb_1,\cdots,\zb_m)$ to $(s_{\zb_1},\cdots,s_{\zb_m})$ is $B_m$-equivariant. In particular, the product of the corresponding reflections remains the same:
\[
s_{\zb_1}\cdots s_{\zb_m}=s_{\zb_1}\cdots s_{\zb_{i-1}}s_{\zb_{i+1}'}s_{\zb_i}s_{\zb_{i+2}}\cdots s_{\zb_m}.
\]

\begin{lem}\label{lem braid}
Let $(\zb_1,\cdots,\zb_m)$ be a sequence of positive roots. Then, for any $i=1,2\cdots,m$ there is a $\zs\in B_m$ so that the first entry in $\zs(\zb_1,\cdots,\zb_m)$ is $\zb_i$.
\end{lem}

\begin{pf}
The element $\zs=(\zs_{i-1}\cdots\zs_1)^{-1}\in B_m$ has the desired property.
\qed\end{pf}
%


\begin{lem}\label{lem group equality implies braid action}
Suppose that $s_{\zb_1}\cdots s_{\zb_m}=s_{\zb_1}\cdots s_{\zb_{p-1}}s_\zg s_{\zb_p}\cdots s_{\zb_{m-1}}$. Then 
\[
\zs_p\zs_{p+1}\cdots \zs_{m-1}(\zb_1,\cdots,\zb_m)=(\zb_1,\cdots,\zb_{p-1},|\zg|,\zb_p,\cdots,\zb_{m-1}).
\]
\end{lem}

\begin{pf}
If we get $\zg'$ instead of $|\zg|$ in the braid equation, then $s_\zg=s_{\zg'}$ which implies that $\zg'$ equals $  \zg$ or $-\zg$ . Since $\zg'\in\Phi_+$ this implies $\zg'=|\zg|$.
\qed\end{pf}

\end{subsection}

\begin{subsection}{Transitive action}
%
We are ready to state the main result of this section.

\begin{thm}\label{thm 2}
Let $W$ be a Coxeter group generated by the simple reflections $s_1,\cdots,s_n$ and let $t_1,\cdots,t_m$ with $m\le n$, be reflections (conjugates of the simple reflections) whose (inverse) product is
\[
	t_mt_{m-1}\cdots t_1=C=s_1s_2\cdots s_n.
\]
Then $m=n$ and there is an element of the braid group $B_n$ which transforms the word $t_n\cdots t_2t_1$ to $s_1s_2\cdots s_n$. That is, the braid group acts transitively on the set of all sequences of $n$ reflections whose product is the Coxeter element $C$.
\end{thm}
%

The proof of the theorem  will be given in the following two subsections.

\end{subsection}

\begin{subsection}{Projective roots}

Take a sequence of $m\le n$ positive roots
\[
	\zb_\ast=(\zb_1,\cdots,\zb_m)
\]
with the property that the product of the corresponding reflections is the inverse Coxeter element:
\[
	s_{\zb_1}\cdots s_{\zb_m}=C^{-1}=s_n\cdots s_1.
\]
Assume that $m$ is minimal. We want to show that $m=n$ and that the action of the braid group as described in section \ref{sect 1.1} is transitive on the set of such sequences. Taking $m$ to be minimal implies that the $\zb_i$ are distinct. Otherwise we could cancel a pair of reflections in the sequence, using the  action of the braid group, and make $m$ smaller.

A priori, $m$ might be smaller than $n$, and the braid group $B_m$ acts on the set of all such sequences $\zb_\ast$ by conjugating the corresponding reflections with each other. The first step is to collect all the projective roots on the left so that, in particular, $\zb_m$ will not be projective, unless all the roots are projective.

Recall that the \emph{projective roots} $p_1,\cdots,p_n$ are the $n$ roots which are sent to negative roots by $C^{-1}$. More precisely, for $i=1,2,\ldots,n$,
\[
	p_i=s_1\cdots s_{i-1}(\za_i).
\]
Let $i$ be minimal so that the projective root $p_i$ occurs in some sequence in the orbit of $\zb_\ast$ under the action of the braid group. By Lemma \ref{lem braid}, we can move $p_i$ to the front to get
\[
	\zb_\ast\sim (p_i,\ast,\cdots,\ast)
\]
where $\sim$ means lying in the same orbit under the action of $B_m$. Let $j>i$ be minimal so that
\[
	\zb_\ast\sim(p_i,p_j,\ast,\cdots,\ast).
\]
Continuing in this way we get the following.

\begin{lem}\label{projective sequence} There is a sequence of positive integers
\[
	j_1<j_2<\cdots<j_k\le n
\]
with $0\le k\le m\le n$
so that
\[
	\zb_\ast\sim(p_{j_1},p_{j_2},\cdots,p_{j_k},\zg_1,\cdots,\zg_{m-k})
\]
satisfying the following conditions.
\begin{enumerate}
\item The action of the braid group $B_{m-k}$ on the sequence of roots $\zg_1,\cdots,\zg_{m-k}$ produces no projective roots. In particular none of the $\zg_i$ is projective.
\item The action of the braid group $B_{m-i}$ on the sequence of roots $p_{j_{i+1}},\cdots,p_{j_k}$, $\zg_1$,$\cdots$,$\zg_{m-k}$ does not produce any projective root $p_j$ with $j<j_{i+1}$.
\end{enumerate}
\end{lem}

\begin{cor}\label{lem 29}
If $k=m$ then $m=n$ and
\[
	\zb_\ast\sim(p_1,p_2,\ldots,p_n)\sim(\za_n,\ldots,\za_2,\za_1).
\]
\end{cor}

\begin{pf}
If $k=m$ then all of the roots are projective. Moreover, since $p_{j_i}=s_1s_2\cdots s_{j_i-1}(\za_{j_i})$, we have $s_{p_{j_i}}s_1s_2\cdots s_{j_i-1} = s_1s_2\cdots s_{j_i-1} s_{j_i}$, by Lemma \ref{lem 1}. Consequently,
\begin{equation}\label{eq 41}
	C=s_1s_2\cdots s_n={s_{p_{j_m}}\cdots s_{p_{j_1}}}s_1\cdots \widehat{s_{j_1}}\cdots \widehat{s_{j_m}} \cdots s_n,
\end{equation}
and since we have $C={s_{p_{j_m}}\cdots s_{p_{j_1}}}$, this implies that
\[
	s_1\cdots \widehat{s_{j_1}}\cdots \widehat{s_{j_m}} \cdots s_n=e\in W
\]
where the notation $\widehat s_i$ means that the reflection $s_i$ is deleted from the sequence. But a product of distinct simple reflections cannot be trivial. (If the product were trivial then the letters could be cancelled two at a time, and the last two letters to be cancelled would be equal.) Therefore, the trivial word is empty and $m=n$ as claimed.
Moreover, the two sides of equation (\ref{eq 41}) lie in the same orbit under the braid group action, which shows that $(p_1,p_2,\ldots,p_n)\sim(\za_n,\ldots,\za_2,\za_1).$
\qed\end{pf}

To prove Theorem \ref{thm 2} it therefore suffices to show that $k=m$.
\end{subsection}

\begin{subsection}{Proof of Theorem \ref{thm 2}}

Let $\zb_\ast=(\zb_1,\ldots,\zb_m)$ be such that $s_{\zb_1}\cdots s_{\zb_m}=C^{-1}=s_ns_{n-1}\cdots s_1$. By Lemma \ref{projective sequence}, we may suppose without loss of generality that
$\zb_\ast= (p_{j_1},p_{j_2},\ldots,p_{j_k}$, $\zb_{k+1},\ldots,\zb_m)$, where $p_{j_1},\ldots,p_{j_k}$ are projective and the two conditions of Lemma \ref{projective sequence} are satisfied.

We want to show that $k=m$.

\underline{Construction}: Define $g_i=s_1s_2\cdots s_{j_i-1}$ for $i=1,2,\ldots,k$. Then $g_i(\za_{j_i})=p_{j_i}=\zb_i$. We can also write $g_i= w_1s_{j_1}w_2s_{j_2}\cdots w_i$, where
 $w_i\in W$ is given by $w_i= s_{j_{i-1}+1}s_{j_{i-1}+2}\cdots s_{j_{i}-1}$ if $i>1$, and $w_1=g_1$. We want to describe the $\zb_i$ with $i>k$ in a similar way. The roots are the orbits of the simple roots under the action of the Coxeter group $W$. Thus, also for $i=k+1,k+2,\ldots,m$, there is a simple root $\za_{j_i}$ and an element $g_i\in W$ such that 
\begin{equation}\label{sign equation}
g_i(\za_{j_i})=\zb_i \textup{ or } g_i(\za_{j_i})=-\zb_i.
\end{equation}
 Let $s_{j_i}$ be the simple reflection corresponding to the simple root $\za_{j_i}$ and let 
 $w_i\in W$ be such that for each $i=1,2,\ldots,m$,
 \begin{equation}\label{eq 31}
g_i=w_1s_{j_1}w_2s_{j_2}\cdots w_{i-1}s_{j_{i-1}}w_i .
\end{equation}
One example to keep in mind is the preprojective case where this sequence of simple reflections can be taken to be the first part of a power of the Coxeter element. Recall that $\beta$ is called \emph{preprojective} if there is some non-negative integer $q$ such that $C^q \beta$ is projective.
  
As consequences of the above construction we have the recursive equations
 \begin{eqnarray}
  g_i&=&g_{i-1}s_{j_{i-1}}w_i \label{eq 331},\\
  s_{\zb_i}g_i &=& g_i s_{j_i}\label{eq 332}, \qquad\qquad \textup{by Lemma \ref{lem 1}}
\end{eqnarray}
  and thus, for $ 2\le i\le m,$ we have
\begin{equation}\label{eq 333}
g_i=s_{\zb_{i-1}}g_{i-1}w_i.
\end{equation}
Now equation (\ref{eq 332}) implies  $g_m s_{j_m} = s_{\zb_m}g_m$, and
applying equation (\ref{eq 333}) repeatedly yields
\[\begin{array}{rcl}
 g_ms_{j_m}&=& s_{\zb_m}s_{\zb_{m-1}}g_{m-1}w_{m} \\
 &=& s_{\zb_m}s_{\zb_{m-1}}s_{\zb_{m-2}}  g_{m-2}w_{m-1}w_{m} \\
  &=& s_{\zb_m}s_{\zb_{m-1}}\cdots s_{\zb_{1}}  g_{1} w_2 \cdots w_{m-1}w_{m} \\
  &=& s_{\zb_m}s_{\zb_{m-1}}\cdots s_{\zb_{1}}  w_1 w_2 \cdots w_{m-1}w_{m} ,
\end{array}\]
  and therefore
\begin{equation}\label{eq 334}
g_ms_{j_m} = C w_1w_2\cdots w_m.
\end{equation}

\underline{Induction hypothesis}: Let $L=(\ell(w_1),\ell(w_2),\ldots,\ell(w_m))$, where $\ell(w_i)$ denotes the length of $w_i$. Now consider the orbit of $(\zb_{k+1},\zb_{k+2},\ldots,\zb_{m})$ under the action of the braid group $B_{m-k}$. We suppose without loss of generality that among all sequences in this orbit our sequence $(\zb_{k+1},\zb_{k+2},\ldots,\zb_{m})$ and choice of $g_i$ for $i>k$ in the construction above are such that the corresponding length vector $L=(\ell(w_1),\ell(w_2),\ldots,\ell(w_m))$ is minimal in lexicographic order.

\begin{lem}\label{lem positive signs}
The signs in equation (\ref{sign equation}) above are all positive, that is, $g_i(\za_{j_i})=\zb_i\in\Phi_+$ for all $1\le i\le m$.
\end{lem}

\begin{pf}
Suppose that $g_i(\za_{j_i})=-\zb_i$. Then, by the Exchange Condition, the expression (\ref{eq 31}) for $g_i$ can be factored as $g_i=as_jb$ where $b(\za_{j_i})=\za_j$ and $a(\za_j)=\zb_i$. Since the $\zb$'s are all distinct, the letter $s_j$ cannot be equal to any of the letters $s_{j_p}$ in (\ref{eq 31}). Therefore, $s_j$ occurs in the middle of some $w_p$ for $p\le i$. So, $w_p=a_ps_jb_p$ with $\ell(a_p)<\ell(w_p)$ and $\zb_i=g_p'(\za_j)$ where
\[
	g_p'=w_1s_{j_1}w_2s_{j_2}\cdots w_{p-1}s_{j_{p-1}}a_p.
\]
But then $\zb_\ast\sim (\zb_1,\dots,\zb_{p-1},\zb_i,\ast,\dots,\ast)$ with corresponding length vector
\[
	L'=(\ell(w_1),\dots,\ell(w_{p-1)},\ell(a_p),\ast,\dots,\ast)<L
\]
contradicting the minimality of $L$.
\qed\end{pf}

\begin{lem}\label{lem 32}
If $m>k$ then $w_1w_2\cdots w_m(\za_{j_m}) $ is a negative root.
\end{lem}

\begin{pf}
We have
\[w_1w_2\cdots w_m(\za_{j_m}) = C^{-1}g_ms_{j_m}(\za_{j_m})=-C^{-1} (\zb_m),\]
where the first equation follows from equation (\ref{eq 334}), and the second follows because  $s_{j_m}(\za_{j_m})=-\za_{j_m}$ and $g_m(\za_{j_m})=\zb_m$ by the previous lemma.

Thus $w_1w_2\cdots w_m(\za_{j_m})$ is a negative root if and only if $C^{-1}\zb_m$ is a positive root, that is, if and only if $\zb_m$ is not projective. This holds, since $m>k$.
\qed
\end{pf}

\begin{lem}\label{lem 33}
If $m>k$ then there exist $p$ with $k<p\le m$ and $a_p, s_i, b_p \in W$ such that $w_p=a_p s_i b_p$ with $\ell(w_p)=\ell(a_p)+\ell(b_p)+1$ and 
\[s_{\zb_1}s_{\zb_2} \cdots s_{\zb_{p-1}}s_{\zg}s_{\zb_{p}}\cdots s_{\zb_{m-1}}=C^{-1},\]
and thus $\zb_\ast\sim (\zb_1,\cdots,\zb_{p-1},|\zg|,\zb_p,\cdots,\zb_{m-1})$ where $\zg = g_{p-1}s_{j_{p-1}}a_p(\za_i)$.
\end{lem}
\begin{pf}
By Lemma \ref{lem 32}, $w_1w_2\cdots w_m$ maps the positive root $\za_{j_m}$ to a negative root, thus there is at least one letter $s_i$ in any expression for $w_1w_2\cdots w_m$ such that $w_1w_2\cdots w_m=as_i b$ and 
$b(\za_{j_m})$ is a positive root and $s_ib(\za_{j_m})$  is a negative root. It follows that $b(\za_{j_m})=\za_i$. We choose as our expression a product of reduced expressions for each $w_i$.
Let $p$ be such that the letter $s_i$ lies in the chosen reduced expression for $w_p$. Then $w_p=a_ps_i b_p$ where $\ell(a_p)<\ell(w_p)$ and
\[w_1w_2\cdots w_m=w_1w_2\cdots w_{p-1} a_ps_ib_pw_{p+1}\cdots w_m.\]
Applying Lemma \ref{lem 1} to the equation $b(\za_{j_m})=\za_i$ yields $s_i b =bs_{j_m}$, which implies
\begin{equation}\label{eq 35}
w_1w_2\cdots w_m = a b s_{j_m}.
\end{equation}

On the other hand,  equation (\ref{eq 334}) yields
$C w_1w_2\cdots w_m = g_ms_{j_m}$, and using equation (\ref{eq 333}) repeatedly, we get
\[C w_1w_2\cdots w_m = s_{\zb_{m-1}}s_{\zb_{m-2}}\cdots s_{\zb_p} g_p w_{p+1}w_{p+2}\cdots w_ms_{j_m}.\]
Now $g_p=g_{p-1}s_{j_{p-1}}w_p = g_{p-1}s_{j_{p-1}} a_ps_ib_p ,$ where the first equality follows from equation (\ref{eq 331}). Since $\zg=g_{p-1}s_{j_{p-1}}a_p(\za_i)$, Lemma \ref{lem 1} implies that $g_{p-1}s_{j_{p-1}}a_ps_i=s_\zg  g_{p-1}s_{j_{p-1}}a_p$, and thus
\[C w_1w_2\cdots w_m = s_{\zb_{m-1}}s_{\zb_{m-2}}\cdots s_{\zb_p}s_\zg g_{p-1}s_{j_{p-1}}a_p b_p w_{p+1}w_{p+2}\cdots w_ms_{j_m}.\]
Again, applying  equation (\ref{eq 333}) repeatedly and using $w_1w_2\cdots w_{p-1}a_pb_p w_{p+1}\cdots w_m = ab$, we get
\[C w_1w_2\cdots w_m = s_{\zb_{m-1}}s_{\zb_{m-2}}\cdots s_{\zb_{p}} s_\zg s_{\zb_{p-1}} \cdots s_{\zb_1} abs_{j_m}.\]
Comparing this result with equation (\ref{eq 35}), we conclude that
\[C=  s_{\zb_{m-1}}s_{\zb_{m-2}}\cdots s_{\zb_{p}} s_\zg s_{\zb_{p-1}} \cdots s_{\zb_1}.\]
Finally, we note that $\zg $ is not projective, by Lemma \ref{projective sequence}, thus, $p>k$.
\qed
\end{pf}

{\it Proof of Theorem \ref{thm 2}}.
Recall that $\zb_\ast=(\zb_1,\zb_2,\ldots,\zb_m)$ was such that $\zb_1,\ldots,\zb_k$ are projective and $\zb_{k+1},\zb_{k+2},\ldots,\zb_m$ are such that 
\[L=(\ell (w_1),\ell(w_2),\ldots,\ell(w_m))\]
is minimal in lexicographic order.

By Corollary \ref{lem 29}, it suffices to show that $m=k$. Suppose $m>k$, then Lemma \ref{lem 33} implies that 
$C=  s_{\zb_{m-1}}s_{\zb_{m-2}}\cdots s_{\zb_{p}} s_\zg s_{\zb_{p-1}} \cdots s_{\zb_1}$ with $\zg= g_{p-1}a_p(\za_i)$. Moreover, 
$s_{\zb_{m-1}}s_{\zb_{m-2}}\cdots s_{\zb_{p}} s_\zg s_{\zb_{p-1}} \cdots s_{\zb_1}$ 
is obtained from the word
\[
w_1s_{j_1}\cdots {a_ps_ib_p}s_{j_p}w_{p+1}\cdots w_{m-1}s_{j_{m-1}}
\]
by pulling out the letters $s_{j_1},\cdots,s_{j_{p-1}},s_i,s_{j_p},\cdots,s_{j_{m-1}}$, note that ${a_ps_ib_p}=w_p$.

 Since $\ell(a_p)<\ell(w_p)$, we conclude that the length vector 
\[L'=(\ell (w_1),\ell(w_2),\ldots,\ell(w_{p-1}),\ell(a_p), \ldots)\]
is strictly smaller than $L$ which contradicts the minimality of $L$. This completes the proof. \qed
\end{subsection}
\end{section}



\begin{section}{Finite type}\label{sect main}

\def\qQM{q^Q(\underline\dim M)}

Let $A$ be an hereditary algebra which is finite dimensional over some field $k$ and of finite representation type, and let $n$ be the number of isoclasses of simple $A$-modules. Then it is well-known that the indecomposable $A$-modules have the same dimension vectors as the representations of an associated modulated quiver $Q$ whose underlying graph is a Dynkin diagram. (See \cite{IT1}, \cite{IT2}.)
Denote the vertices of 
$Q$ by $1,2,\ldots,n$, and let $W$ be the corresponding Weyl group. 
Note that $W$ is a finite crystallographic reflection group, and all finite crystallographic reflection groups appear in that way. 
Thus $W$ is the Coxeter group generated by $S=\{s_1,s_2,\ldots,s_n\}$ subject to the relations $(s_is_j)^{m_{ij}}=1$, where $m_{ii}=1$, and for $i\ne j$, we have $m_{ij}=2 $ if there is no edge $i- j$ in the Dynkin diagram, and $m_{ij}=3,4,6$ if the weight of the edge $i-j$ in the Dynkin diagram is $1,2,3$ respectively.

The Weyl group $W$ acts on the root space $\mathbb R^n$ of $Q$ and this action is conjugate to the action of $W$ on $\mathbb R^n$ defined in the last section. The details, which we do not need, are as follows. There is a $W$-equivariant linear isomorphism $\varphi:\mathbb R^n\to \mathbb R^n$ from the root space of $Q$ to the root space of $W$ given by
\[
	\varphi(x_1,\cdots,x_n)=\left(\sqrt{f_1}x_1,\cdots,\sqrt{f_n}x_n\right)
\]
where $f_i=\dim_k S_i=\dim_k\End_A(S_i)$ with $S_i$ the $i$th simple $A$-module. This induces a $W$-equivariant bijection from the root system $\Phi^Q$ of $Q$ to the root system $\Phi^W$ of $W$ by sending the root $\zb\in\Phi^Q$ to $\varphi(\zb)$ divided by its length.

By  \cite{G}, \cite{DR}, the dimension vector $\underline{\textup{dim}}$ gives a 1-1 correspondence between isomorphism classes of indecomposable $A$-modules $M$ and positive roots of $Q$ sending the isomorphism class $[M]$ of $M$ to 
\[
\underline{\textup{dim}}\, M=(d_1,\ldots,d_n)=d_1\za_1+\cdots+d_n\za_n.
\]
The corresponding positive root of $W$ is $\varphi(\underline\dim M)$ divided by its length which is $\sqrt{\dim_k\End_A(M)}$. However, the results of the last section are not affected by rescaling the positive roots or changing the basis for the root space. In this section we will reinterpret these results using the root system of $Q$.

\begin{subsection}{Adapted expressions} 

A sequence of reflections $s_{i_1}s_{i_2}\cdots s_{i_m}$  is called \emph{adapted} to the quiver $Q$ if $i_1$ is a sink of $Q$, and $i_k$
is a sink of the quiver $s_{i_{k-1}}\cdots s_{i_2}s_{i_1}Q$, for each $k\ge 2$, where $s_i Q$ is the quiver obtained from $Q$ by reversing all arrows at vertex $i$.

Since $Q$ has no cycles, we can assume without loss of generality that the Coxeter element $C=s_1s_2\cdots s_n$ is adapted to the quiver $Q$. Let $\nu$ be the number of positive roots, and let $s_{i_1}s_{i_2}\cdots s_{i_\nu}$ be an adapted reduced expression of the longest element $w_0$ of the Weyl group.
It is well known, see for example \cite[Ch IV, 1.6]{Bou}, that the sequence 
\begin{equation}\label{eq 42}
(\za_{i_1},s_{i_1}(\za_{i_2}),\ldots ,s_{i_1}s_{i_2}\cdots s_{i_{\nu-1}}(\za_{i_\nu}) )
\end{equation}
 contains every positive root exactly once. This induces a total order on the positive roots by $\za<\zb$ if $\za $ appears in the sequence (\ref{eq 42}) before $\zb$.

\begin{rem}
If $M$ is a non-projective indecomposable $A$-module, then the dimension vector $\underline{\textup{dim}}\, \tau M$ of its Auslander-Reiten translate $\tau M$ is equal to $C \underline{\textup {dim}} \, M$.
\end{rem}

\begin{rem}\label{rem homorder} Let $M,N$ be indecomposable $A$-modules such that
 $\underline{\dim}\,M< \underline{\dim}\,N$. Then
$\Hom_A(N,M)=0$ and $\Ext_A(M,N)=0$.
\end{rem}

\begin{pf}
Indeed, suppose to the contrary that $\Hom_A(N,M)\ne 0$. Since in finite representation type every module is preprojective, we can suppose without loss of generality that $M$ is projective (otherwise apply the Coxeter transformation repeatedly to $M$ and $N$ until the image of $M $ is projective). Since  $\Hom_A(N,M)\ne 0$ and $A$ is hereditary, it follows that $N$ is projective too and that $N$ is a submodule of $M$. Let  $\za^s=\underline{\dim}\,M, \za^t=\underline{\dim}\,N$, so by hypothesis we have $s<t$. Then $\za^s=s_{i_1}s_{i_2}\cdots s_{i_{s-1}} \za_{i_s}$ and these $s-1$ simple reflections are without repetition and similarly $\za^t=s_{i_1}s_{i_2}\cdots s_{i_{t-1}} \za_{i_t}$ and these $t-1$ simple reflections are also without repetition. Consequently, the support of $M$ is a subset of the support of $N$, and since both are projective, it follows that $M$ is a submodule of $N$. Thus $M=N$, a contradiction.
The vanishing of $Ext_A(M,N)$ follows from the Auslander-Reiten formula.\qed
\end{pf}

A special case is $M=\tau N$, where $\tau$ is the Auslander-Reiten translation. Then the Auslander-Reiten formula yields $\Ext_A(N,N)=D\Hom_A(N,M)=0$, reflecting the fact that	 all indecomposable $A$-modules are exceptional.
\end{subsection}

\begin{subsection}{Exceptional sequences}\label{sect es}
We recall some facts about exceptional sequences. These results in this subsection are also valid if the hereditary algebras $A$ is not of finite representation type.
\begin{definition}\label{def es}
 A sequence of modules $(E_1,\cdots,E_r)$ is an \emph{exceptional sequence} if
\begin{enumerate}
\item $\End_A(E_i)$ is a division algebra for all $i$,
\item $\Hom_A(E_j,E_i)=0$ for $j>i$,
\item $\Ext_A(E_j,E_i)=0$ for $j\ge i$.
\end{enumerate}
If $r=n$ then $(E_1,\cdots,E_n)$ is called a \emph{complete exceptional sequence}. When $r=1$, $E_1$ is called an \emph{exceptional module}.
\end{definition}

For example, the projective modules $(P_1,\cdots,P_n)$ form a complete exceptional sequence and the simple modules in reverse order $(S_n,\cdots,S_1)$ form a complete exceptional sequence. Since $A$ is hereditary, we know by \cite{HR} that, for any indecomposable module $E$ with $\Ext_A(E,E)=0$, $\End_A(E)$ must be a division algebra. 
Therefore, condition (1) may be replaced with the assumption that each $E_i$ is indecomposable.

The braid group acts on the set of complete exceptional sequences as follows.
 The generator $\zs_i$ of the braid group (which moves the $i$-th strand over the $i+1$st strand) acts on a complete exceptional sequence  $E=(E_1,\cdots,E_n)$ by
\[
	\zs_iE=(E_1,\cdots,E_{i-1},X,E_i,E_{i+2},\cdots,E_n)
\]
where $X$ is the unique module making the indicated sequence exceptional. 
See \cite{C-B93} for details. Note that our action is the inverse of \cite{C-B93} since we prefer the label of the strand that goes under to change. 
Moreover, the dimension vector of $X$ is given by
\begin{equation}\label{eq 1}
\underline{\textup{dim}}\,X= \,s_{e_i}(e_{i+1}), \textup{ or } 
\underline{\textup{dim}}\,X= \,-s_{e_i}(e_{i+1}),
\end{equation}
where $e_i,e_{i+1}$ are the positive roots corresponding to the dimension vectors of the modules $E_i, E_{i+1}$ respectively.  

\begin{thm}[Crawley-Boevey \cite{C-B93}, Ringel \cite{R94}] \label{thm CB} The braid group acts transitively on the set of (isomorphism classes of) exceptional sequences. \qed\end{thm}

\begin{cor}\label{cor CB}
Let  $s_{e_1}s_{e_2}\cdots s_{e_n} \in W$, denoted $s_E$, be the product of the reflections corresponding to the dimension vectors of the elements of an exceptional sequence $E$. 
Then for any two complete exceptional sequences $E,E'$ we have $s_E=s_{E'}$.
\end{cor}
\begin{pf}
Applying Lemma \ref{lem 1} to equation (\ref{eq 1}), we get
\[
	s_{\underline{\textup{dim}}\,X}=s_{e_i}s_{e_{i+1}}s_{e_i}.
\]
Thus, for each generator $\zs_i$ of the braid group, we have $s_{\zs_i E}=s_{ E}$, since  $$\cdots s_{e_{i-1}} (s_{e_{i}} s_{e_{i+1}} s_{e_i}) s_{e_{i}} s_{e_{i+2}} \cdots
= \cdots s_{e_{i-1}} s_{e_{i}} s_{e_{i+1}} s_{e_{i+2}} \cdots $$
The statement now follows from the Theorem.
\qed
\end{pf}
\end{subsection}

\begin{subsection}{Cluster categories}
%
We suppose without loss of generality that the last $n$ positive roots in the sequence (\ref{eq 42}) are the dimension vectors of the indecomposable injective $A$-modules in order from $1$ to $n$, that is, for $k=1,2\ldots,n$, we have $\underline{\textup{dim}} \,I(k)= 
s_{i_1}s_{i_2}\cdots s_{i_{\nu-n+k-1}}(\za_{i_{\nu-n+k}})$.

Let $(j_1,j_2,\ldots, j_{n+\nu})$ be the sequence $(i_1,\ldots,i_\nu)$ with $(1,\ldots,n)$ inserted at the beginning, thus 
\[ (j_1,j_2,\ldots, j_{n})=(1,2,\ldots, n) \quad \textup{and} \quad  (j_{n+1},j_{n+2},\ldots, j_{n+\nu})= (i_1,i_2,\ldots,i_\nu).\]
Let $w_1=Cw_0=s_{j_1}s_{j_2}\cdots s_{j_{n+\nu}}$. Note that this is an adapted expression which is not reduced. For $t=1,2,\ldots,\nu+n$, define
\[ \za^t=s_{j_1}s_{j_2}\cdots s_{j_{t-1}}(\za_{j_t}).
\]

Then $\za^1,\za^2,\ldots,\za^\nu$ are precisely the positive roots and 
$\underline{\textup{dim}} \,P(k)=\za^{k}$ if $1\le k\le n$.
Moreover,
 $\za^{\nu+1},\za^{\nu+2},\ldots,\za^{\nu+n}$ are the negative roots $\za^{\nu+k}=-\underline{\textup{dim}} \,P(k)$ for $1\le k\le n$.

The cluster category $\CA$, introduced in \cite{BMRRT}, \cite{CCS}, is the orbit category of the derived category $\DA$ under the endofunctor $\tau^{-1}[1]$, where $\tau$ denotes the Auslander-Reiten translation and $[1]$ is the shift.  As a fundamental domain for $\CA$, we may take $\ind A \cup A[1]$, in other words, every indecomposable object in $\CA$  is the orbit of an indecomposable $A$-module or of the first shift of an indecomposable projective $A$-module, see \cite{BMRRT}.
  
 For $t=1,2,\ldots,\nu$, let $M_t$ be the indecomposable $A$-module whose dimension vector is equal to $\za^t$, and for $t=\nu+k$, with $k>0$, let $M_t=P(k)[1]$ be the first shift of the indecomposable projective $A$-module $P(k)$. Then the indecomposable objects in $\mathcal{C}_A$ are in bijection with $M_1,M_2,\ldots,M_{\nu+n}$.

Let $N_t$ denote the indecomposable $A$-module $N_t=M_t$ if $t\le \nu $, and $N_t=P(k)$ if $t = \nu+k$ with $k>0$.
Then the dimension vector of $N_t$ is $\za^{t}$ if $t\le \nu$, and it is $\za^{k}=-\za^t$ if $t=\nu+k$, with $k>0$.
In particular,
\begin{equation}\label{eq 2}
s_{\underline{\textup{dim}\,}N_{t}}=s_{\za^{t}} \quad \textup{ if  $t\le \nu$}.
\end{equation}

Let  $w^\zd(t_1,\ldots,t_n)$ denote the expression that is obtained from the expression $s_{j_1}\cdots s_{j_{n+\nu}}$ for $w_1$ by deleting the reflections at the positions $t_1,t_2,\ldots,t_n$.

\begin{thm}\label{thm main}
Let $1\le t_1<t_2<\cdots <t_n\le \nu+n$. Then the following statements are equivalent. 
\begin{enumerate}
\item $M=M_{t_1}\oplus M_{t_2}\oplus\cdots\oplus M_{t_n}$ is a cluster-tilting object in $\mathcal{C}_A$.
\item $N=(N_{t_1},N_{t_2},\ldots,N_{t_n})$ is an exceptional sequence in $\textup{mod}\,A$.
\item $s_{\za^{t_n}}s_{\za^{t_{n-1}}}\cdots s_{\za^{t_1}}=
s_1s_2\cdots s_n$.
\item $w^\zd(t_1,\ldots,t_n)$ is a reduced expression for $w_0$.
\end{enumerate}
\end{thm}

\begin{rem}
Note that the theorem implies that there  is a map from cluster-tilting objects to exceptional sequences, which is obtained by ordering the indecomposable summands of the cluster-tilting object and then replacing the shifts of projectives (if any) by the corresponding projectives.
This map is neither injective nor surjective! 

To see that it is not injective, it suffices to take one cluster-tilting object $T=P(1)\oplus\cdots\oplus P(n) $ to be the sum of the indecomposable projective modules and another cluster-tilting object $T'=T[1]$ to be the sum of the shifts of the indecomposable projective modules. Both are mapped to the same exceptional sequence $(P(1),\ldots,P(n))$. 

To see that the map is not surjective, it suffices to notice that the exceptional sequence $(S_n,\ldots,S_1)$ consisting of the simple modules in reverse order is not in the image. 

In the theorem, we need to fix a sequence of integers $t_i$ to have the equivalence of {\rm (1)} and {\rm (2)}.
\end{rem}
The proof of Theorem \ref{thm main} will be given in section \ref{sect proof}.
\end{subsection}

\begin{subsection}{Mutations}

In this subsection, we give a precise description of the mutations in terms of the reduced expressions. This will be seen to be more or less equivalent to Theorem \ref{thm main}.

\begin{lem}\label{rem 1}
If $w^\zd(t_1,t_2,\ldots,t_n)$ is a reduced expression for $w_0$, then for each $k$ there is a unique $t_k' \in\{1,2,\ldots,n+\nu\}\setminus \{t_1,t_2,\ldots,t_n\}$ such that
\[w^\zd(t_1,\ldots,t_{k-1},t_{k+1},\ldots,t_n,t_k') \]
 is a reduced expression for $w_0$. 
\end{lem} 

We will say that $t_k'$ is obtained from $t_k$ by \emph{algebraic mutation}.

\begin{pf} The claim is that, if the letter $t_k$ is inserted into its original place in the deleted word $w^\zd(t_1,t_2,\ldots,t_n)$, then there is a unique other letter $t_k'$ which needs to be deleted in order for the word to remain equal to $w_0$ in the Coxeter group.
This is a special case of the following statement. If we insert a simple reflection $t_k$ in the middle of any reduced expression for the longest word $w_0$ in any finite Coxeter group, say $w_0=ab\mapsto at_kb$, then there exists a unique other letter which needs to be removed in order for the result to remain equal to $w_0$. The letter that needs to be removed, call it $t_k'$, is in either $a$ or $b$ but not both so that either
\begin{enumerate}
\item $a=a_1t_k' a_2$ and $w_0=a_1a_2t_kb$ or
\item $b=b_1t_k' b_2$ and $w_0=at_kb_1b_2$.
\end{enumerate}
To see this, suppose that $t_k=s_i$ and let $\za=a(\za_i)$ and $\zb=b^{-1}(\za_i)$. Then $\za=w_0(\zb)$. Since $w_0$ sends all positive roots to negative roots and vice versa, exactly one of the roots $\za,\zb$ is positive and the other is negative. If $\za$ is negative, then $at_k$ is not reduced and, by the exchange condition, there is a unique letter $t_k'$ in $a$ so that $a=a_1t_k'a_2$ and $at_k=a_1a_2$, so we are in case 1 above. If $\zb$ is negative we are in case 2.
\qed\end{pf}

Define 
\[\overline{s}_{j_i} =\left\{ \begin{array}{ll}
1 &\textup{if $i \in\{t_1,t_2,\ldots,t_n\} $}\\
s_{j_i} & \textup{otherwise.}
\end{array}\right.\]
Then $w^\zd(t_1,t_2,\ldots,t_n)=\overline{s}_{j_1}\overline{s}_{j_2}\cdots
\overline{s}_{j_{n+\nu}}$. 

Let $k\in\{1,2,\ldots,n\}$. Consider two cases.
\begin{enumerate}
\item  Suppose first that there exists a positive integer  $\ell < t_k$ such that 
\[\overline{s}_{j_\ell}\overline{s}_{j_{\ell+1}}\cdots
\overline{s}_{j_{t_k -1}} (\za_{j_{t_k}}) \] 
is a negative root, and let $\ell$  be the largest such integer. Then 
\[\overline{s}_{j_{\ell+1}}\cdots
\overline{s}_{j_{t_k -1}} (\za_{j_{t_k}})=\za_\ell,\] and, by Lemma \ref{lem 1}, 
$\overline{s}_{j_{\ell+1}}\cdots
\overline{s}_{j_{t_k -1}}\overline{s}_{j_{t_k }} 
=
\overline{s}_{j_{\ell}}\overline{s}_{j_{\ell+1}}\cdots
\overline{s}_{j_{t_k -1}}$. 
Consequently
\[w^\zd(t_1,\ldots,\ell,\ldots t_{k-1},t_{k+1},\ldots,t_n) =w^\zd(t_1,\ldots t_{k-1},t_k,t_{k+1},\ldots,t_n).
\]
Thus $w^\zd(t_1,\ldots,\ell,\ldots t_{k-1},t_{k+1},\ldots,t_n) $ is a reduced expression of $w_0$ and, hence, $\ell $ is the unique $t_k'$ in Lemma \ref{rem 1}.

\item Now suppose that there exists a positive integer $\ell  > t_k$ such that 
\[\overline{s}_{j_\ell}\overline{s}_{j_{\ell-1}}\cdots
\overline{s}_{j_{t_k +1}} (\za_{j_{t_k}}) \] 
is a negative root, and let $\ell$  be the least such integer.
Then 
\[\overline{s}_{j_{\ell-1}}\cdots
\overline{s}_{j_{t_k +1}} (\za_{j_{t_k}})=\za_\ell,\] and, again by Lemma \ref{lem 1}, 
$\overline{s}_{j_{\ell-1}}\cdots
\overline{s}_{j_{t_k +1}}\overline{s}_{j_{t_k }} 
=
\overline{s}_{j_{\ell}}\overline{s}_{j_{\ell-1}}\cdots
\overline{s}_{j_{t_k +1}}$. 
Consequently, 
$\overline{s}_{j_{t_k +1}}\cdots\overline{s}_{j_{\ell-1}}\overline{s}_{j_{\ell}}= \overline{s}_{j_{t_k }} \overline{s}_{j_{t_k +1}}\cdots\overline{s}_{j_{\ell-1}}$, and
\[w^\zd(t_1,\ldots, t_{k-1},t_{k+1},\ldots,\ell,\ldots,t_n) =w^\zd(t_1,\ldots t_{k-1},t_k,t_{k+1},\ldots,t_n).
\]
Thus $w^\zd(t_1,\ldots t_{k-1},t_{k+1},\ldots,\ell,\ldots,t_n) $ is a reduced expression of $w_0$ and, hence, $\ell $ is the unique $t_k'$ in Lemma \ref{rem 1}.
\end{enumerate}
It also follows from Lemma \ref{rem 1} that exactly one of the two cases above must hold.
Assuming Theorem \ref{thm main} we obtain the following theorem.

\begin{thm}\label{thm 3}

Let $T=M_{t_1}\oplus\cdots\oplus M_{t_n}$ be a tilting object. Let $\mu_k(T)=(\oplus_{j\ne k}M_{t_j})\oplus M_{t_k'}$ be the mutation of $T$ in direction $k$. Then either $t_k'<t_k$ is  the largest integer  such that 
\[\overline{s}_{j_\ell}\overline{s}_{j_{\ell+1}}\cdots
\overline{s}_{j_{t_k -1}} (\za_{j_{t_k}}) \] 
is a negative root, or  $t_k'>t_k$ is  the least integer such that 
\[\overline{s}_{j_\ell}\overline{s}_{j_{\ell-1}}\cdots
\overline{s}_{j_{t_k +1}} (\za_{j_{t_k}}) \] 
is a negative root.
\end{thm}

\end{subsection}

\begin{subsection}{An example}\label{example}
Let $Q$ be the quiver $ 1\ot 2\ot 3\ot 4$ of type $A_4$, let $c=s_1s_2s_3s_4$ and use the reduced expression $s_1s_2s_1s_3s_2s_1s_4s_3s_2s_1$ for $w_0$. Then $w_1$ is the word
\[\begin{array}{ccccccccccccccccccccc}
s_1&s_2&s_3&s_4&s_1&s_2&s_1&s_3&s_2&s_1&s_4&s_3&s_2&s_1\\
1&2&3&4&5&6&7&8&9&10&11&12&13&14
\end{array}\]
where the second row indicates the position for convenience.
Let $T$ be the tilting object in the cluster category whose direct summands are the indecomposable projective modules, that is, $(t_1,t_2,t_3,t_4)=(1,2,3,4)$ and $T=M_1\oplus M_2 \oplus M_3 \oplus M_4$. The corresponding reduced expression $w^\zd(1,2,3,4)$ equals
\[\begin{array}{ccccccccccccccccccccc}
&&&&s_1&s_2&s_1&s_3&s_2&s_1&s_4&s_3&s_2&s_1\\
1&2&3&4&5&6&7&8&9&10&11&12&13&14
\end{array}\]
We can mutate $T$ in the four positions $1,2,3$ and $4$, and the resulting reduced expressions in order are 

\[\begin{array}{ccccccccccccccccccccc}
s_1&&&&&s_2&s_1&s_3&s_2&s_1&s_4&s_3&s_2&s_1\\
&s_2&&&s_1&s_2&&s_3&s_2&s_1&s_4&s_3&s_2&s_1\\
&&s_3&&s_1&s_2&s_1&s_3&s_2&&s_4&s_3&s_2&s_1\\
&&&s_4&s_1&s_2&s_1&s_3&s_2&s_1&s_4&s_3&s_2&\\
1&2&3&4&5&6&7&8&9&10&11&12&13&14
\end{array}\]

If we mutate the second row in position $3$, we obtain the following reduced expression corresponding to $\mu_{t_3}\mu_{t_2} T$:
\[\begin{array}{ccccccccccccccccccccc}
&s_2&s_3&&s_1&s_2&&s_3&&s_1&s_4&s_3&s_2&s_1\\
1&2&3&4&5&6&7&8&9&10&11&12&13&14
\end{array}\]

These reduced expressions together with their positions describe the combinatorial structure of (a fundamental domain of) the Auslander-Reiten quiver of the corresponding cluster-tilted algebras, compare with \cite{B}. The last reduced expression, for example, gives rise to the quiver
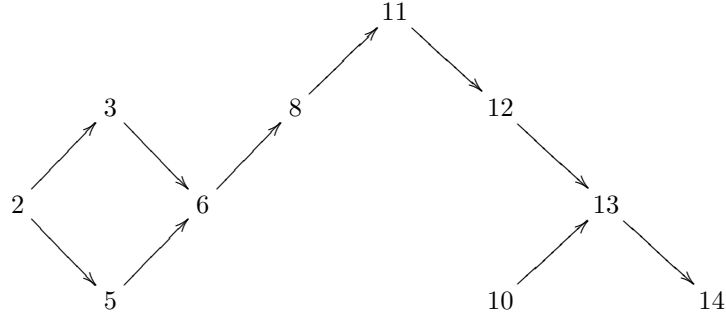
\[\xymatrix{
&&&&11 \ar[rd]\\
& 3 \ar[rd] && 8 \ar[ru]&&12 \ar[rd]\\
2 \ar[ru] \ar[rd]&& 6 \ar[ru]&&&&13  \ar[rd]\\
&5 \ar[ru] && &&10 \ar[ru]&&14\\
}
\]
Where the labels are the positions of the ten simple reflections $s_i$ in the reduced expression, and label $k$ is placed at the level $i$ in the quiver if the reflection at position $k$ is $s_i$. For example, the reflection  $ s_2 $ appears at positions $2,6,13 $ in the reduced expression and therefore the labels 
$2,6,13 $ form the second level of  the Auslander-Reiten quiver.
To obtain the complete Auslander-Reiten quiver, one has to add one arrow  $13 \to 2$ and one arrow $14\to 3$. 
This can  be obtained by forming the long (not reduced) word 
\begin{equation}\label{long word} w^\zd(t_1,t_2,\ldots,t_n) 
\rho(w^\zd(t_1,t_2,\ldots,t_n)) w^\zd(t_1,t_2,\ldots,t_n) \cdots\end{equation} 
and identifying vertices with labels $i$ and $j$ whenever $i\equiv j \,(\textup{mod} \, n+\nu)$,
where $\rho$ is the endomorphism of the Weyl group given by
$ \rho(s_i)= s_{\tilde\rho(i)}$, and $\tilde \rho$ is the identity in types $A_1, D_{n}$, with $n$ odd, $ E_7,$ and $E_8$;  and $\tilde \rho$ is the unique non-trivial automorphism of the Dynkin diagram in the types $A_n$ with $n>1$, $D_{n}$ with $n$ even and $E_6$. Thus in our example the long word is 
\[\begin{array}{cccccccccccccccccccccccccccccccccccccccccccccc}
&s_2&s_3&&s_1&s_2&&s_3&&s_1&s_4&s_3&s_2&s_1 \\1&2&3&4&5&6&7&8&9&10&11&12&13&14\\
\\
&s_3 & s_2&&s_4&s_3&&s_2&&s_4 &s_1&s_2&s_3&s_4&\cdots
\\15&16&17&18&19&20&21&22&23&24&25&26&27&28&\cdots
\end{array}\]

In this example, the reflections that correspond to the indecomposable projective modules can be found in the reduced expression $w^\zd(1,4,7,9)$ as follows. The deleted positions $1,4,7,9$ correspond to simple reflections $s_{j_1}=s_1,s_{j_4}=s_4,s_{j_7}=s_1,s_{j_9}=s_2$ respectively. To each of these positions $i$, 
let $p_i>i$ be the smallest integer such that $\overline{s}_{p_i}=s_{j_i}$ (in the long word in  (\ref{long word})). Thus 
$p_1=5, p_4=11, p_7=10, p_9 =13.$ The  indecomposable projective modules are at the position $(p_1,p_4,p_7,p_9)=(5,11,10,13)$ in the Auslander-Reiten quiver.

\end{subsection}

\begin{subsection}{A geometric interpretation}\label{Woo}

The $A_n$ example gives a cluster-tilting theoretic interpretation of the following types of drawings which are described in \cite{Woo04}. The $14$ crossings in Figure \ref{fig1} correspond to the simple reflections 
\[
	(s_1),s_2,s_1,s_3,s_2,(s_4),(s_1),s_3,(s_2),s_4,s_1,s_3,s_2,s_1.
\]
 Since the Weyl group of $A_4$ is the symmetric group $S_5$, the simple reflections are the simple transpositions $s_i=(i,i+1)$ and the product of these simple transpositions is the permutation indicated in Figure \ref{fig1}. If the simple transpositions in parentheses are deleted then we get the longest word which is the permutation indicated in Figure \ref{fig2}.

The mutation process is easy to visualize, we simply take two lines in Figure \ref{fig2} and make them cross where they do not and make them not cross where they do.

\def\crossx{\put(-.5,-.5){\line(1,1){1}}\put(.5,-.5){\line(-1,1){1}}}
\def\crossplus{\qbezier(-.5,-.5)(0,0)(.5,.5)
	\put(.5,-.5){\line(-1,1){.4}}\put(-.5,.5){\line(1,-1){.4}}}
\def\crossplusright{\put(-.1,-.1){\line(1,1){.6}}\put(.1,-.1){\line(1,-1){.4}}}
\def\crossminus{\qbezier(.5,-.5)(0,0)(-.5,.5)
	\put(.5,.5){\line(-1,-1){.4}}\put(-.5,-.5){\line(1,1){.4}}}
\def\crossnot{\crosscup\crosscap}
\def\crosscap{\qbezier(-.5,-.5)(0,0)(.5,-.5)}
\def\crosscup{\qbezier(-.5,.5)(0,0)(.5,.5)}
\newcommand{\halfcapright}[1]{\qbezier(-.1,.1)(-.5,.5)(-1,.5)\put(-1,.5){\line(-1,0){#1}}}
\newcommand{\halfcupright}[1]{\qbezier(-.1,-.1)(-.5,-.5)(-1,-.5)\put(-1,-.5){\line(-1,0){#1}}}
\newcommand{\halfcapleft}[1] {\qbezier(.1,.1)(.5,.5)(1,.5)\put(1,.5){\line(1,0){#1}}}
\newcommand{\halfcupleft}[1] {\qbezier(.1,-.1)(.5,-.5)(1,-.5)\put(1,-.5){\line(1,0){#1}}}

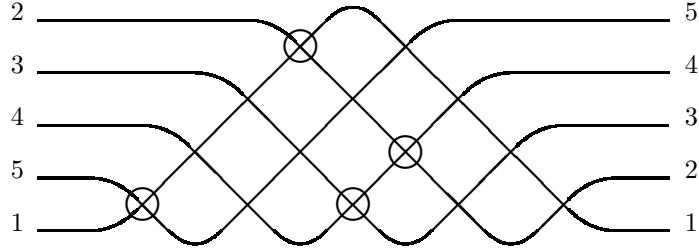
\begin{figure}[htbp]
\begin{center}
{
\setlength{\unitlength}{.7cm}
{\mbox{
\begin{picture}(10,5)
      \thicklines
    \put(1,1){\halfcupright{1}\halfcapright{1}\put(-.1,.1){\line(1,-1){.6}}\put(-.1,-.1){\line(1,1){3.6}}}
      \put(2,2){\halfcapright{2}\put(-.1,.1){\line(1,-1){.6}}}
      \put(3,3){\halfcapright{3}\put(-.1,.1){\line(1,-1){.6}}}
      \put(4,4){\halfcapright{4}\put(-.1,.1){\line(1,-1){.6}}}
   \put(9,1){\halfcupleft{1}\put(.1,-.1){\line(-1,1){3.6}}}
   \put(9,1){\halfcapleft{1}\put(.1,.1){\line(-1,-1){.6}}}
      \put(7,1){\put(.5,-.5){\line(-1,1){.7}}}
      \put(8,2){\halfcapleft{2}\put(.1,.1){\line(-1,-1){.7}}}
            \put(7,3){\halfcapleft{3}\put(.1,.1){\line(-1,-1){.6}}}
                        \put(6,4){\halfcapleft{4}\put(.1,.1){\line(-1,-1){.6}}}
                             \put(3,1)\crossx
    \put(5,1){\crossx}
    \put(5,3){\crossx}
    \put(6,2){\crossx}
    \put(7,1){\crossx}
    \put(1,1){\circle{.6}}
    \put(5,1){\circle{.6}}
    \put(6,2){\circle{.6}}
    \put(4,4){\circle{.6}}
    \put(2,0)\crosscup
    \put(4,0)\crosscup
   \put(6,0)\crosscup
   \put(8,0)\crosscup
    \put(4,2)\crossx
 	\put(5,5)\crosscap
	\put(-1.5,.5){1}
	\put(-1.5,1.5){5}
	\put(-1.5,2.5){4}
	\put(-1.5,3.5){3}
	\put(-1.5,4.5){2}
	\put(11.3,.5){1}
	\put(11.3,1.5){2}
	\put(11.3,2.5){3}
	\put(11.3,3.5){4}
	\put(11.3,4.5){5}
\end{picture}}
}}
\caption{The positions corresponding to $p_1=5, p_4=11, p_7=10, p_9 =13$ in the previous example are circled.}
\label{fig1}
\end{center}
\end{figure}
\begin{figure}[htbp]
\begin{center}
%
{
\setlength{\unitlength}{.6cm}
{\mbox{
\begin{picture}(10,5)
      \thicklines
    \put(1,1){\qbezier(-.5,-.5)(-.8,-.5)(-2,-.5)\qbezier(-.5,.5)(-.8,.5)(-2,.5)\put(.5,.5){\line(1,1){2}}}
    \put(1,1){\crossnot}
      \put(2,2){\halfcapright{2}\put(-.1,.1){\line(1,-1){.6}}}
      \put(3,3){\halfcapright{3}\put(-.1,.1){\line(1,-1){.6}}}
      \put(4,4){\qbezier(-.5,.5)(-1,.5)(-5,.5)}
     \put(4,4){\crossnot}
   \put(9,1){\halfcupleft{1}\put(.1,-.1){\line(-1,1){3.6}}}
   \put(9,1){\halfcapleft{1}\put(.1,.1){\line(-1,-1){.6}}}
      \put(7,1){\put(.5,-.5){\line(-1,1){.7}}}
      \put(8,2){\halfcapleft{2}\put(.1,.1){\line(-1,-1){.7}}}
      \put(7,3){\halfcapleft{3}\put(.1,.1){\line(-1,-1){.6}}}
            \put(6,4){\halfcapleft{4}\put(.1,.1){\line(-1,-1){.6}}}
     \put(3,1)\crossx
    \put(5,1){\crossnot}
    \put(5,3){\crossx}
    \put(6,2){\crossnot}
    \put(7,1){\crossx}
    \put(2,0)\crosscup
    \put(4,0)\crosscup
   \put(6,0)\crosscup
   \put(8,0)\crosscup
    \put(4,2)\crossx
 	\put(5,5)\crosscap
	\put(-1.5,.5){5}
	\put(-1.5,1.5){4}
	\put(-1.5,2.5){3}
	\put(-1.5,3.5){2}
	\put(-1.5,4.5){1}
	\put(11.3,.5){1}
	\put(11.3,1.5){2}
	\put(11.3,2.5){3}
	\put(11.3,3.5){4}
	\put(11.3,4.5){5}
\end{picture}}
}}
\caption{Every pair of lines crosses exactly once, making this a reduced expression for the longest word $w_0=(15)(24)\in S_5$.}
\label{fig2}
\end{center}
\end{figure}
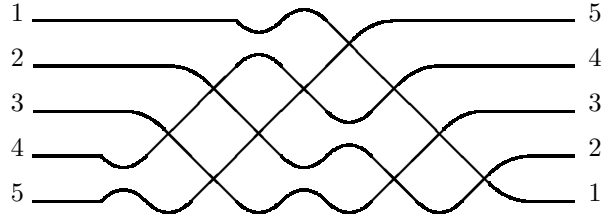

\end{subsection}

\end{section}


\begin{section}{Proof of Theorem \ref{thm main}}\label{sect proof}


The proof of the Theorem is subdivided into the three Lemmas \ref{lem 12} -- \ref{lem 23}. 
\begin{subsection}{Cluster-tilting objects and exceptional sequences}\label{sect 12}

\begin{lem}\label{lem 12}
Let $1\le t_1<t_2<\cdots <t_n\le \nu+n$. Then $M=M_{t_1}\oplus M_{t_2}\oplus\cdots\oplus M_{t_n}$ is a cluster-tilting object in $\mathcal{C}_A$ if and only if  $N=(N_{t_1},N_{t_2},\ldots,N_{t_n})$ is an exceptional sequence in $\textup{mod}\,A$.
\end{lem}

\begin{pf}
Suppose $M=M_{t_1}\oplus\cdots\oplus M_{t_n}$ is a tilting object in $\CA$ and let $t_i<t_j$. 
Recall that $N_t=M_t$ if $1\le t\le \nu$, and if $\nu< t \le \nu +n$ then $M_t$ is the object $P(t-\nu)[1]=\tau P(t-\nu)$ in the cluster category and $N_t$ is the module $P(t-\nu)$.
If $t_i<t_j\le\nu$, then $\Hom_A(N_{t_j},N_{t_i})=0$ by, Remark \ref{rem homorder}.
If $t_i\le\nu<t_j$, then $\Hom_A(N_{t_j},N_{t_i})=\Hom_A(P(t_j-\nu),M_{t_i})$. If the latter is nonzero then so is
\[\Hom_{\CA}(P(t_j-\nu),M_{t_i})\cong D\Ext_{\CA}(M_{t_i}, \tau P(t_j-\nu)) = D\Ext_{\CA}(M_{t_i}, M_{t_j}),\]
which is impossible, since $M$ is a tilting object. Thus, again, we have $\Hom_A(N_{t_j},N_{t_i})=0$.
Finally, if $\nu< t_i<t_j$, then $\Hom_A(N_{t_j},N_{t_i})=\Hom_A(P(t_j-\nu),P(t_i-\nu))$, which is zero by Remark \ref{rem homorder}.
Therefore, we have $\Hom_A(N_{t_j},N_{t_i})=0$ if $t_j>t_i$.

If $t_j\le\nu$ then 
$\Ext_A(N_{t_j},N_{t_i})=0 $ because $\Ext_{\CA}(M_{t_j},M_{t_i})=0$, and if $\nu<t_j$ then 
$\Ext_A(N_{t_j},N_{t_i})=\Ext_A(P(t_j-\nu),N_{t_i})=0$, since $P(t_j-\nu)$ is projective. 
Thus $(N_{t_1},N_{t_2},\ldots,N_{t_n})$ is an exceptional sequence.

Conversely, suppose that $(N_{t_1},N_{t_2},\ldots,N_{t_n})$ is an exceptional sequence.
Suppose first that $t_i<t_j\le \nu$. Then 
\[D\Ext_{\CA}(N_{t_j},N_{t_i}) \cong \Ext_{\CA}(N_{t_i},N_{t_j})
\cong \Ext_A(N_{t_i},N_{t_j})\oplus \Ext_A(N_{t_j},N_{t_i}),\] 
where the second summand is zero because  $(N_{t_1},N_{t_2},\ldots,N_{t_n})$ is an exceptional sequence, and the first summand is isomorphic to
$D\Hom_A(N_{t_j},\tau N_{t_i})$ which is zero by Remark \ref{rem homorder}.
Thus $\Ext_{\CA}(M_{t_j},M_{t_i}) = \Ext_{\CA}(M_{t_i},M_{t_j}) = 0$. 

Now suppose $t_i \le\nu< t_j$. Then 
\[\begin{array}{rcl}
\Ext_{\CA}(M_{t_j},M_{t_i}) &=&\Ext_{\CA}(P(\nu-t_j)[1],N_{t_i})\\
&\cong& \Hom_{\CA}(P(\nu-t_j), N_{t_i})\\
&=&  \Hom_{\DA}(N_{t_j}, N_{t_i})\oplus \Hom_{\DA}(P(\nu-t_j), \tau^{-1} N_{t_i}[1]),\\
\end{array}\]
where the first summand is zero because $N$ is an exceptional sequence and the second summand is zero because of the structure of the derived category.

Finally, suppose that $\nu<t_i<t_j$, then 
\[
\Ext_{\CA}(M_{t_j},M_{t_i}) =\Ext_{\CA}(P(\nu-t_j),P(\nu-t_i)) = 0.\\
\]
Thus $M=\oplus_{i=1}^n N_{t_i}$ has no self-extension in $\CA$, whence $M$ is a tilting object in $\CA$. 
This completes the proof.
\qed\end{pf}
\end{subsection}

\begin{subsection}{Reduced expressions and the Coxeter element}\label{sect 34}

\begin{lem}\label{lem 34}
Let $1\le t_1<t_2<\cdots <t_n\le \nu+n$. Then  $w^\zd(t_1,\ldots,t_n)$ is a reduced expression for $w_0$ if and only if
$s_{\za^{t_n}}s_{\za^{t_{n-1}}}\cdots s_{\za^{t_1}}=
s_1s_2\cdots s_n$.
\end{lem}
\begin{pf}
Setting $w=s_{j_1}s_{j_2}\cdots s_{j_{t-1}}$
and $\za=\za_{j_t}$, we have $w(\za)=\za^t$, and then Lemma \ref{lem 1} implies
\[ s_{\za^t}= s_{j_1}s_{j_2}\cdots s_{j_{t-1}}s_{j_t}s_{j_{t-1}} \cdots s_{j_2}s_{i_1}.\]
Therefore $s_{\za^{t_1}}s_{\za^{t_{2}}}\cdots s_{\za^{t_n}}$ is equal to
\[\begin{array}{cl}
& s_{j_1}s_{j_2}\cdots s_{j_{t_1-1}} \widehat{s_{j_{t_1}}} s_{j_{t_1+1}} \cdots s_{j_{t_2-1}} \widehat{s_{j_{t_2}} } s_{j_{t_2+1}} \cdots  s_{j_{t_n-1}} s_{j_{t_n}}  s_{j_{t_n-1}} s_{j_{t_n-2}} \cdots s_{j_1},
\end{array}\]
and multiplying with $w=s_{j_1}s_{j_2}\cdots s_{j_{\nu+n}}$ on the right, we get
\[(s_{\za^{t_1}}s_{\za^{t_{2}}}\cdots s_{\za^{t_n}})\, w_1= w^\zd(t_1,\ldots,t_n).\]
Since $w_1=s_1 s_2\cdots s_n\,w_0$, it follows that $w^\zd(t_1,\ldots,t_n)$ is a reduced expression for $w_0$ if and only if 
  $s_1s_2\cdots s_n =s_{\za^{t_n}}s_{\za^{t_{n-1}}}\cdots s_{\za^{t_1}}$, and this completes the proof.
\qed\end{pf}

\end{subsection}

\begin{subsection}{Exceptional sequences and the Coxeter element}\label{sect 23}

\begin{lem} \label{lem 23}
Let $1\le t_1<t_2<\cdots <t_n\le \nu+n$. Then
 $N=(N_{t_1},N_{t_2},\ldots,N_{t_n})$ is an exceptional sequence in $\textup{mod}\,A$ if and only if 
$s_{\za^{t_n}}s_{\za^{t_{n-1}}}\cdots s_{\za^{t_1}}=
s_1s_2\cdots s_n$.
\end{lem}
\begin{pf}
Let $E=(S_n,S_{n-1},\ldots,S_1)$ be the exceptional sequence given by the simple modules in reverse order.
By Theorem \ref{thm CB}, there is an element $\zs$ in the braid group such that $\zs(E)=N$ is the exceptional sequence under consideration, and Corollary \ref{cor CB} implies 
\[s_ns_{n-1}\ldots s_1=s_E = s_{N} = s_{\za^{t_1}} s_{\za^{t_2}}\cdots s_{\za^{t_n}},\]
where the last identity follows from equation (\ref{eq 2}).

Conversely, $s_E=s_ns_{n-1}\cdots s_1 =C^{-1}$ and by our assumption this is equal to $s_{\za^{t_1}}s_{\za^{t_2}}\cdots s_{\za^{t_n}}$. Then Theorem \ref{thm 2} yields the existence of an element of the braid group $\zs$ such that $\zs s_E= s_{\za^{t_1}}s_{\za^{t_2}}\cdots s_{\za^{t_n}}$, and, hence,  $\zs E= (N_{t_1},N_{t_2},\ldots,N_{t_n}) $ is an exceptional sequence by Theorem \ref{thm CB}. 
\qed
\end{pf}
\end{subsection}

\begin{subsection}{Alternate proof}\label{subsec 24}

We note that Theorem \ref{thm 2} was not used in the proofs of Lemma \ref{lem 12} and Lemma \ref{lem 34} and was only used in the second half of Lemma \ref{lem 23} above. In the finite case, this can be replaced by Lemma \ref{rem 1} in the following way. 

Let $(t_1,\cdots,t_n)$ be minimal in lexicographic order so that $s_{\za^{t_1}}s_{\za^{t_2}}\cdots s_{\za^{t_n}}=C^{-1}$ but we do not know if the corresponding objects form a cluster. Then we can use Lemma \ref{rem 1} to algebraically mutate the last term $t_n$ to $t_n'$ so that
$
	w^\zd(t_1,\ldots,t_{n-1},t_n') =w_0
$
and therefore
\[
	s_{\za^{t_1}}s_{\za^{t_2}}\cdots s_{\za^{t_k}}s_{\za^{t_n'}}s_{\za^{t_{k+1}}} \cdots s_{\za^{t_{n-1}}}=C^{-1}
\]
Since $t_n'\neq t_n$ we must have that $k<n-1$. Therefore $(t_1,\cdots,t_k,t_n',t_{k+1},\cdots,t_{n-1})$ is less than $(t_1,\cdots,t_n)$ in lexicographic order. Thus, by induction, the objects $M_{t_1},\cdots,M_{t_{n-1}},M_{t_n'}$ form a cluster. The object $M_{t_n'}$ can be mutated to an object, say $M_{t_n^\ast}$, to obtain another cluster $M_{t_1},\cdots,M_{t_{n-1}},M_{t_n^\ast}$ where $t_n^\ast\neq t_n'$.

We claim that $t_n^\ast=t_n$ proving the lemma and thus the theorem. By Theorem \ref{thm CB} we know that the product of reflections corresponding to this new cluster is $C^{-1}$. By Lemma \ref{lem 34} this implies that $w^\zd(t_1,\ldots,t_{n-1},t_n^\ast) =w_0$. But this equation determines $t_n^\ast$ uniquely by Lemma \ref{rem 1}. So, $t_n^\ast=t_n$ as claimed.

\end{subsection}

\end{section}


\begin{section}{Infinite type}\label{sect infinite type}


We will now look at the analogue of Theorem \ref{thm main} for quivers of infinite type. For simplicity of terminology we restrict to the simply laced case. 

\begin{subsection}{Quivers of infinite type}\label{sect 41}

Suppose that $Q$ is a quiver without oriented cycles and $K$ is an algebraically closed field. The path algebra $KQ$ is a finite dimensional hereditary algebra over $K$. Kac \cite{K80} showed that the dimension vectors of the indecomposable $KQ$-modules are exactly the positive roots of the Kac-Moody Lie algebra associated to $KQ$. The Weyl group $W$ of $KQ$ is generated by reflections with respect to the bilinear form $B$ given by $B(\za_i,\za_i)=1$ and $B(\za_i,\za_j)=-n_{ij}/2$ where $n_{ij}=n_{ji}$ is the number of arrows between vertices $i$ and $j$. This is a Coxeter group (See, e.g. \cite{Bro}, $[A]\Rightarrow[C]$). So, Theorem \ref{thm 2} applies. The relations among the simple reflections $s_i$ are $(s_is_j)^{m_{ij}}=1$ where $m_{ij}=2,3$ for $n_{ij}=0,1$ respectively and $m_{ij}=\infty$ when $n_{ij}\ge2$. The bilinear form $B$ for $Q$ agrees with the bilinear form defined for $W$ if and only if all of the numbers $n_{ij}$ are $\le2$.

The \emph{Coxeter element} $C\in W$ is given by the product of simple reflections
\[
	C=s_1s_2\cdots s_n
\]
which we assume as before to be adapted to the quiver $Q$.

Kac \cite{K80} defines a \emph{real root} of $Q$ to be any root of the form $w(\za_i)$ where $w\in W$ and $\za_i$ is a simple root and he showed that a root $\za$ is real if and only if $B(\za,\za)=1$. If all $n_{ij}\le2$ then the set of real roots of $Q$ is equal to the set of roots of $W$. If there are any $n_{ij}\ge3$ then the root spaces of $Q$ and $W$ will not be $W$-equivariantly isomorphic. However, there is still a $W$-equivariant bijection between the set of positive roots of $W$ and the set of positive real roots of $Q$ since both sets are in $W$-equivariant bijection with the set $T$ of section \ref{sect 1.1} with a positive real root $\zb$ corresponding to the reflection $s_\zb\in T$ given as before by
\[
	s_\zb(x)=x-2B(\zb,x)\zb.
\]

\end{subsection}
\begin{subsection}{Exceptional sequences}\label{sect 42}

Recall from Kac \cite{K82} and Schofield \cite{S92} that a \emph{real Schur root} is a real root which is also the dimension vector of an indecomposable $KQ$-module $M$ so that $\End_{KQ}(M)=K$. Since $\dim M$ is a real root, this implies that $M$ is an exceptional module. Conversely, the dimension vectors of all exceptional modules are real Schur roots.

Combining Theorem \ref{thm CB} of Crawley-Boevey and our Theorem \ref{thm 2} we get the following.

\begin{thm}\label{thm 41}
Suppose that $\zb_1,\zb_2,\ldots,\zb_n$ is a sequence of real roots of $W$. Then the following are equivalent.
\begin{enumerate}
\item There is an exceptional sequence $(E_1,E_2,\cdots,E_n)$ with $\underline\dim E_i=\zb_i$.
\item The product of the corresponding reflections is the inverse of the Coxeter element: $s_{\zb_1}s_{\zb_2}\cdots s_{\zb_n}=C^{-1}$.
\end{enumerate}
\end{thm}

This has the following corollary where we recall that a \emph{prefix} of the Coxeter element is defined to be any element $w\in W$ which can be expressed as a product of reflections $w=t_1t_2\cdots t_k$ for which there exist $n-k$ other reflections $t_{k+1}, \cdots,t_n$ so that
\[
	t_1t_2\cdots t_n=C.
\]

\begin{cor}
A real root $\zb$ is a real Schur root if and only if $s_{\zb}$ is a prefix of the Coxeter element.
\end{cor}

\begin{pf}
If $\zb$ is a real Schur root, then $\zb=\underline\dim E$ for some exceptional module $E$. This can be completed to an exceptional sequence $(E_1,\cdots,E_n)$ with $E_n=E$ by \cite{C-B93}. Therefore $C=s_\zb s_{\zb_{n-1}}\cdots s_{\zb_1}$ by Theorem \ref{thm 41} where $\zb_i=\underline\dim E_i$, so $s_\zb$ is a prefix. 

Conversely, if $s_\zb$ is a prefix, say $s_\zb=t_1\cdots t_k$ then $C=s_\zb t_{k+1}\cdots t_n$. By \cite{Dyer}, the Coxeter element cannot be written as a product of fewer than $n$ reflections,  which implies that $k=1$. So, Theorem \ref{thm 41}  implies that $\zb$ is the dimension vector of an exceptional module, so $\zb$ is a real Schur root.\qed
\end{pf}

\end{subsection}

\section*{Acknowledgements} The authors would like to thank G. Todorov, H. Thomas, W. Crawley-Boevey, A. Hubery and R. Marsh for very stimulating helpful discussions, as well as the referees for helping to improve the exposition. This work was developed right before and right after the Auslander Memorial Lectures and Conference at Woods Hole with some crucial suggestions by Gordana Todorov and the results were improved during ICRA 13 at S\~{a}o Paulo, Brazil with substantial help from Hugh Thomas.

\end{section}

\section*{Appendix (with Hugh Thomas)}
Hugh Thomas has communicated to us the following application of our work to his joint paper with Ingalls \cite{Ing-Th}. 

 Let $Q$ be a quiver. In \cite{Ing-Th}, the authors define a map $\phi$ from finitely generated, exact abelian, extension-closed subcategories of $\textup{mod}\, kQ$ to prefixes of the Coxeter element. The purpose of this appendix is to show that this map is a bijection in general type, generalizing the result of \cite{Ing-Th} where it has been shown that it is a bijection in finite and affine type.

The map $\phi$ is defined by taking an exceptional sequence in the subcategory and then considering the product of the corresponding reflections. Since the exceptional sequence in the subcategory can be extended to one in the full category, the product of reflections is a prefix of the Coxeter element.

\begin{thm}\label{thmhugh}
The map $\phi $ is a bijection from finitely generated, exact abelian, extension-closed subcategories of $\textup{mod}\, kQ$ to prefixes of the Coxeter element.
\end{thm}

\begin{pf}
First we show surjectivity.  Let $w=t_1\dots t_{r}$ be a prefix of $c$.  So there is some $u=t_{r+1}\dots t_n$  such that $wu=c$.  By Theorem \ref{thm 41} above, the reflections $t_i$ are reflections in real Schur roots, and the corresponding objects form an exceptional sequence.  In particular, the minimal exact abelian and extension-closed subcategory containing the objects corresponding to the reflections $t_1,\dots,t_r$ is sent to $w$ by the map $\phi$.

Next we show injectivity.  Suppose there were two subcategories $A,B$ which map to the same prefix of $c$, say $w=t_1\dots t_r$. Then $t_1,\dots,t_r$ can be extended to a factorization of $c$ into $n$ reflections, say, by adding reflections $t_{r+1},\dots,t_n$.  By  Theorem \ref{thm 41} above, this implies that the exceptional sequences for $A$ and $B$ can be extended to an exceptional sequence for the full category, by appending the objects corresponding to $t_{r+1},\dots,t_n$.  But this means that $A$ and $B$ are both the perpendicular to the subcategory generated by the objects corresponding to the reflections $t_{r+1},\dots,t_n$.  Thus $A=B$.
\qed
\end{pf}

\end{document}